\documentclass{amsart}
\usepackage{amssymb}
\usepackage{verbatim}
\usepackage{array}
\usepackage{latexsym}
\usepackage{enumerate}
\usepackage{booktabs}

\usepackage{amsmath,mathtools}
\usepackage{mathrsfs}
\usepackage{amsfonts}
\usepackage{amsthm}
\usepackage[foot]{amsaddr}
\usepackage[table,svgnames,dvipsnames]{xcolor}
\usepackage{caption}
\usepackage[english]{babel}
\usepackage{graphicx}
\usepackage{tikz}
\usepackage{tikz-cd}
\usepackage{adjustbox}
\usepackage{pgfplots}
\pgfplotsset{compat=newest}
\usepackage{nicematrix}
\usepackage[autostyle]{csquotes}

\usepackage[backend=biber]{biblatex}
\usepackage[autostyle]{csquotes}
\usepackage{hyperref}
\usepackage[capitalise]{cleveref}

\newtheorem{theorem}{Theorem}[section]

\theoremstyle{definition}
\newtheorem{example}[theorem]{Example}
\newtheorem{remark}[theorem]{Remark}
\newtheorem{definition}[theorem]{Definition}

\newcommand{\Fq}{{\mathbb{F}_q}}

\newcommand{\F}{\mathbb{F}}

\newcommand{\N}{\mathbb{N}}

\newcommand{\X}{\mathcal{X}}
\newcommand{\Y}{\mathcal{Y}}
\newcommand{\M}{\mathcal{M}}

\newcommand{\D}{\mathcal{D}}
\newcommand{\sL}{\mathcal{L}}
\newcommand{\Fs}{\mathcal{F}}

\newcommand{\Div}{\mathrm{Div}}
\newcommand{\Supp}{\mathrm{Supp}}
\newcommand{\Span}{\mathrm{Span}}

\newcommand{\RS}{\mathrm{RS}}
\newcommand{\val}{v}
\renewcommand{\div}{\mathrm{div}}

\renewcommand{\P}{\mathbb{P}}

\title{Weierstrass semigroups and the order bound}
\date{}

\author[A. Barraud]{Alix Barraud$^1$} \address{$^1$Institut de Mathématiques de Bordeaux, Université de Bordeaux, Talence, France} \email{alix.barraud@math.u-bordeaux.fr} \thanks{}
\author[Y. Çakıroğlu]{Yağmur Çakıroğlu$^2$} \address{$^2$Department of Mathematics, Hacettepe University, Ankara, Turkey} \email{yagmur.cakiroglu@hacettepe.edu.tr} \thanks{}
\author[B. Gouthier]{Bianca Gouthier$^3$} \address{$^3$Max Planck Institute for Mathematics, Bonn, Germany} \email{gouthier@mpim-bonn.mpg.de} \thanks{}
\author[G.L. Matthews]{Gretchen L. Matthews$^4$} \address{$^4$Department of Mathematics, Virginia Tech,  Blacksburg, VA 24061, USA} \email{gmatthews@vt.edu} \thanks{}
\author[L. Vicino]{Lara Vicino$^{5,6}$} \address{$^5$Faculty of Science and Engineering - Bernoulli Institute, University of Groningen, Groningen 9747AG, The Netherlands}  \email{l.vicino@rug.nl} \thanks{}
\address{$^6$Department of Applied Mathematics and Computer Science, Technical University of Denmark, Kongens Lyngby 2800, Denmark}

\begin{document}

\begin{abstract}
    The aim of this survey is to provide the reader with an essential and accessible introduction to the theory of Weierstrass semigroups, in the context of the theory developed by K.-O. Stöhr and J.F. Voloch. Furthermore, we discuss an application of Stöhr-Voloch theory in coding theory, namely the Feng-Rao bound (also known as the \emph{order bound}) for the dual minimum distance of one-point algebraic geometry codes from a curve, which relies on the knowledge of certain Weierstrass semigroups of the curve. 
\end{abstract}

\maketitle

\thanks{{\em Keywords}: algebraic geometry code, curves over finite fields, Feng-Rao bound, order bound, Weierstrass semigroup.

\thanks{{\em MSC codes}: 11G20, 14G50, 94B27.

\section{Introduction}

Over the last four decades, Weierstrass semigroups of curves over finite fields have been a rich research topic, attracting the attention of a multitude of researchers both in the field of positive-characteristic algebraic geometry and in that of coding theory. In 1986, K.-O. Stöhr and J.F. Voloch published \cite{stohr1986weierstrass}, a milestone paper that played a crucial role in sparking a renewed interest in Weierstrass semigroups within the community of researchers interested in algebraic curves over finite fields. Let $\X$ be a projective, geometrically irreducible, non-singular algebraic curve of genus $g$ defined over the finite field $\Fq$ with $q$ elements, and let $\X(\Fq)$ denote the set of $\Fq$-rational points of $\X$. The celebrated Hasse-Weil bound  
\begin{equation}
\label{eq:hasse-weil}
    ||\X(\Fq)|-(q+1)| \leq 2g\sqrt{q}
\end{equation}
states, in particular, that
\begin{equation}
\label{eq:hasse-weil-upper}
    |\X(\Fq)| \leq q + 1 + 2g\sqrt{q}.
\end{equation}
Curves attaining this upper bound are said to be \emph{maximal}. They are of special interest not only because they are extremal objects with special geometric properties but also for their applications in coding theory. 
While the Hasse-Weil bound is, in general, sharp, under certain assumptions the bound can be improved. For instance, when $q$ is not a square, Serre's bound
\begin{equation*}
    ||\X(\Fq)|-(q+1)| \leq g [2\sqrt{q}],
\end{equation*}
where $[\cdot]$ denotes the integer part, improves on \cref{eq:hasse-weil}. 

In \cite{stohr1986weierstrass}, the authors devise a general framework to tackle the problem of finding improvements to the upper bound in \cref{eq:hasse-weil-upper}. Their idea is to work with explicit coordinates, embedding the curve $\X$ in some projective space and considering the equations of the osculating hyperplanes at the points of $\X$. In this way, they manage to define a function with a controlled number of poles and zeros of high multiplicity at the $\Fq$-rational points of the curve. From this, they derive an upper bound for $|\X(\Fq)|$ that depends on $q$, on the genus $g$, on the dimension of the ambient projective space, on the degree of $\X$, and on the so-called \emph{Weierstrass order-sequence} of the embedding. In particular, by appropriately choosing the embedding, they manage to obtain another proof of the Hasse-Weil bound, and also to obtain improvements on \cref{eq:hasse-weil-upper} in several cases. 

The techniques used in \cite{stohr1986weierstrass} rely on the theory of Weierstrass order-sequences associated to projective embeddings, which is derived by translating into a geometric setting the approach of F.K. Schmidt \cite{schmidt1939arithmetischen}. However, by developing this theory, K.-O. Stöhr and J.F. Voloch also obtain some new results and new insights with respect to \cite{schmidt1939arithmetischen}.

The Weierstrass order-sequence of a curve is strictly connected to the Weierstrass semigroups of the curve, and in fact it can be deduced from the knowledge of the \emph{generic semigroup}, that is the Weiestrass semigroup that all but a finite number of points of the curve have. This motivates why the study of Weierstrass semigroups is of crucial importance within the theory of Weierstrass order-sequences, which also paved the way for further applications in the study of curves over finite fields. For instance, Weierstrass point theory played a key role in the investigation of the uniqueness (up to isomorphism) of maximal curves having the largest or the second-largest genus (see  \cite{AT}, \cite{FGT}, \cite{RS}) which were major open problems. Furthermore, the study of Weierstrass semigroups has also been of importance for the determination of the automorphism groups of curves in several cases; see, for instance, the recent publications \cite{BMV23}, \cite{BMV26}, \cite{beelen2026weierstrass}, \cite{beelen2024some}, \cite{N2025}. In general, understanding the Weierstrass semigroups at the points of a curve embedded in a certain projective space gives not only a rich amount of information on the geometry of the embedded model, but also insights on intrinsic algebraic properties of the curve. 

This motivates why, after the geometric theory of Weierstrass order-sequences was exposed in \cite{stohr1986weierstrass}, the interest in the study of Weierstrass semigroups of curves increased considerably, leading to several papers where the semigroups of curves with special properties (for instance, maximal curves) were explicitly determined. For example, all the Weierstrass semigroups have been determined for the Hermitian curve \cite{GV}, for the Giulietti-Korchmáros curve \cite{BM18GK}, for the Suzuki curve \cite{BMZ21}, for the Skabelund curve obtained as a Kummer cover of the Suzuki curve \cite{BLM21}, for the three known maximal curves having the third-largest genus \cite{BMV23}, \cite{BMV26}, \cite{beelen2026weierstrass}, for some maximal Fermat curves \cite{beelen2026weierstrasssemigroupsmaximalfermat}. The Weierstrass semigroups at certain points of several other curves have been studied as well, see for instance \cite{bartoli2018ag} and \cite{montanucci2020ag}, where Weierstrass semigroups at all the $\F_{q^2}$-rational points of, respectively, the Garcia-Güneri-Stichtenoth curves and the Beelen-Montanucci curves are determined.

Moreover, the applications of curves, and in particular of maximal curves, in coding theory introduced by V.D. Goppa in the series of papers \cite{G70}, \cite{G71}, \cite{G77}, \cite{G81}, \cite{G82}, \cite{G88} turned out to be connected with the study of Weierstrass semigroups too. Indeed, in \cite{GL92}, A. Garcia and R. F. Lax demonstrated how the Weierstrass gap set, meaning complement of the Weierstrass semigroup, may be used to define one-point algebraic geometry codes  (AG codes for short) with minimum distance greater than the usual lower bound. Soon thereafter, the same authors in a collaboration with S. J. Kim \cite{GKL93} proved that consecutive elements of the Weierstrass gap set  
improves the usual lower bound on the minimum
distance of certain one-point AG codes. Furthermore, in \cite{FR94} G.-L. Feng and T.R.N. Rao presented a lower bound on the dual minimum distance of one-point AG codes from curves, which relies on the knowledge of the Weierstrass semigroup at the point contained in the support of the one-point divisor defining the code. The Feng-Rao bound is also known as the \emph{order bound} because of its connection with the theory of order functions (see \cite{HLP98}) and, more specifically, with Weierstrass order-sequences. This bound is always at least as good as the dual Goppa bound and, in fact, it significantly improves on it in several instances. 

Motivated by the study of multi-point AG codes, natural generalisations of Weierstrass semigroups to the case of tuples of points have also been introduced and studied in the literature, see for instance 
\cite{castellanos2024weierstrass},
\cite{HK01}, \cite{HK06}, \cite{HK05}, \cite{M04},\cite{M01}, \cite{M04HE},   \cite{M05},  \cite{s2025set}.
Moreover, a generalised version of the order bound was also derived in \cite{OrderBound}. This bound gives good estimates for the dual minimum distance of multi-point AG codes in several examples, see for instance \cite{TwoPointCodes}, \cite{LTV24}, \cite{LV22}.

The overarching aim of this survey is to provide an essential introduction to the theory of Weierstrass order-sequences developed by K.-O. Stöhr and J.F. Voloch in \cite{stohr1986weierstrass}, which motivates the interest in the study of Weierstrass semigroups of curves over finite fields. To this aim, we highlight key points of the paper by F. Torres \cite{torres2000} which contains a thorough and detailed discussion of the results in \cite{stohr1986weierstrass}. Moreover, we discuss the Feng-Rao bound for the minimum distance of duals of one-point AG codes. These topics provide the reader with an idea not only of the implications that Weierstrass semigroups have in understanding the geometry and the intrinsic properties of curves, but also a glimpse into their applications in coding theory, which are of remarkable interest in the theory of AG codes.

In an ongoing research project, we are working on the computation of Weierstrass semigroups at the $\F_{q^{2n}}$-rational points of the Garcia-Güneri-Stichtenoth curves $GGS_n$ \cite{GGS2010}, for $n\geq 5$. For $n=3$ all the Weierstrass semigroups are already known, see \cite{BMGKcurve2018}, while for $n\geq 5$ only those at the $\F_{q^2}$-rational points of $GGS_n$ are known (see \cite{guneri2013automorphism}, \cite{bartoli2018ag}). 

The survey is organised as follows. In \cref{sec:linearseries}, we introduce Weierstrass semigroups, discussing linear series on a curve and associated morphisms to projective spaces. In this framework, we describe the Hermitian invariants and give the definition of Weierstrass semigroups, pointing out their connections to the geometry of the curve. We include several examples to clarify the exposition.
In \cref{sec:orderbound}, after recalling fundamental definitions on AG codes from curves, we discuss one-point AG codes and the Feng-Rao bound, highlighting the role played by Weierstrass semigroups in the study of the dual minimum distance. The paper concludes with a summary and open problems in Section \ref{section:conclusion}.

\section{Linear series, Hermitian invariants and Weierstrass semigroups}
\label{sec:linearseries}

\subsection{Preliminaries}

We start by fixing notation and recalling some relevant background. Throughout this section, $\mathcal{X}$ will be a projective, irreducible, non-singular algebraic curve over an algebraically closed field $\mathbb{F}$.  Let $\mathbb{F}(\X)$ denote the associated algebraic function field. Moreover, $\Div(\X)$ will denote the group of divisors on $\X$, that is the $\mathbb{Z}$-free abelian group generated by the points of $\X$, whose elements are of the form $D=\sum{n_P}P$. We use the notation $v_P(D):=n_P$ when convenient. We say that a divisor $D$ is \textit{effective} (denoted by $D\succeq0$) if $v_P(D)\geq 0$ for all $P$. Moreover, there is a partial order on $\Div(\X)$: for $D,E\in \Div(\X)$, $D\succeq E$ if and only if $D-E$ is effective.

\begin{remark}
    \label{rem:algclosedfield}
    The assumption that $\F$ is algebraically closed is instrumental for a smoother development of the theory. However, the theory applies also if $\X$ is defined over a finite field $\F_q$. Indeed, in that case one considers the curve as the algebraic curve obtained by extending the constant field to the algebraic closure $\overline{\F}_q$ of $\F_q$, equipped with the action of the Frobenius map. The \emph{rational} points in $\X(\overline{\F}_q)$ are then exactly those fixed by the Frobenius map, and a divisor is \emph{rational} if and only if it is invariant under the Frobenius action.
\end{remark}

To any non-zero rational function $f\in\F(\X)^*$ one can associate the divisor $$\mathrm{div}(f):=\sum v_P(f)P=\div_0(f)-\div_\infty(f)$$ where $v_P$ denotes the valuation at the point $P\in\X$, and $\div_0(f):=\sum\limits_{v_P(f)>0} v_P(f)P$ and $\div_\infty(f):=\sum\limits_{v_P(f)<0} (-v_P(f))P$ are respectively the \textit{zero} and \textit{pole} divisor of $f$. Recall that two divisors $D,E$ are said to be linearly equivalent (denoted by $D\sim E$) if $D=E+\mathrm{div}(f)$ for some $f\in\F(\X)^*$.
Given a divisor $E$ on $\X$, we set $$|E|:=\{D\in\mathrm{Div}(\X)\mid D\succeq0,D \sim E\}=\{E+\mathrm{div}(f)\mid f\in \mathcal{L}(E)\backslash\{0\}\}$$ where $\mathcal{L}(E)$ is the Riemann-Roch space $\mathcal{L}(E):=\{f\in\F(\X)^*\mid E+\mathrm{div}(f)\succeq0\}\cup\{0\}$. Recall that $\mathcal{L}(E)$ is a vector space over $\F$ and that by Riemann-Roch theorem we have that $\ell(E)=\deg(E)+1-g+\ell(C-E)$, where $C$ is a canonical divisor on $\X$, $g$ denotes the genus of $\X$ and $\ell(E):=\dim_\F(\mathcal{L}(E))$.
Notice that $|E|$ can be identified with $\P(\mathcal{L}(E))$ via the map $E+\div(f)\mapsto[f]$. Indeed two rational functions $f,g\in\F(\X)^*$ are such that $\div(f)=\div(g)$ if and only if $f=ag$ for some $a\in\F^*$. We will then keep in mind that $|E|$ has the structure of a projective space and we will use the identification $|E|\cong\P(\mathcal{L}(E))$.
A \textit{linear series} on $\X$ is a subset $$\D:=\{E+\mathrm{div}(f)\mid f\in\D'\backslash\{0\}\}\subseteq|E|$$ where $\D'$ is a linear subspace of $\mathcal{L}(E)$. We will also use the identification $\D\cong\P(\D')$. The \textit{degree} and (projective) \textit{dimension} of a linear series $\D$ are respectively $d=\deg(\D):=\deg(E)=\sum n_P$ and $r=\dim(\D):=\dim_{\F}(\D')-1$. We then say that $\D$ is a $g_d^r$ on $\X$; notice that the dimension of the \textit{empty} linear series is $-1$, since it is obtained for $\D'=\{0\}$. A linear series $\D$ is said to be \textit{complete} if $\D=|E|$.

For $\D$ a linear series, its \textit{multiplicity} at a point $P\in\X$ is $b(P):=\min\{v_P(D)\mid D\in\D\}$. We have that $b(P)>0$ if and only if $P\in\Supp(D)=\{P\in\X\mid v_P(D)\neq0\}$ for all $D\in\D$ (recall that all divisors in a linear series are effective), thus $b(P)\neq0$ for finitely many points $P\in\X$ and we can define the effective divisor $B=B^{\D}:=\sum b(P)P$, called \textit{base locus} of $\D$. The linear series $\D$ is called \textit{base-point-free} if $B=0$. A point in the support of the base locus $B$ is called \textit{base point} of $\D$.

\medskip

\subsection{Hermitian invariants and Weierstrass semigroups}
For $P \in \X(\mathbb{F})$ and any natural number $i$, we can define $$\D_i(P):=\{D\in\D\mid D\succeq iP\}.$$ Then $\D_i(P)$ is a linear series, $\D_i(P)\supseteq\D_{i+1}(P)$ and $\dim(\D_i(P))\leq\dim(\D_{i+1}(P))+1$ (see \cite[Lemma 1.3]{torres2000}). 
 Moreover, $\D_i(P)'=\D'\cap L(E-iP)$. As a consequence, $\D_i(P)=\emptyset$ whenever $i>d$, since in this case $\deg(E-iP)<0$ and thus $L(E-iP)=\{0\}$. Now, if we take $i$ to be the multiplicity $b(P)$, by its definition we have that $\D_{b(P)}(P)=\D$ and, moreover, $\D_{b(P)}(P)\supsetneq D_{b(P)+1}(P)$. 

An integer $b(P)\leq j\leq d$ is called \textit{$(\D,P)$-order} (or \textit{Hermitian $P$-invariant}) if $\D_j(P)\supsetneq D_{j+1}(P)$. In particular, for every $(\D,P)$-order $j$ we have that $$\dim(\D_{j+1}(P))=\dim(\D_j(P))-1,$$ while for integers $i$ that are \textbf{not} $(\D,P)$-orders the dimension does not decrease, and indeed $\D_{i+1}(P)=\D_i(P)$. Therefore, $(\D,P)$-orders detect jumps in the chain $$\D=\D_{b(P)}(P)\supsetneq \D_{b(P)+1}(P)\supseteq\D_{b(P)+1}(P)\supseteq\dots\supseteq\emptyset.$$ Since every strict inclusion makes the dimension drop exactly of one and that we start from dimension $r$ and we get to dimension $-1$, we have exactly $r+1$ $(\D,P)$-orders $$b(P)=j_0<j_1<\dots<j_r.$$ Notice that $j_i$ is the multiplicity at $P$ of the linear series $\D_{j_i}(P)$, for $i=0,\dots,r$, and the latter is a $g_d^{r-i}$ on $\X$.

\begin{remark}\label{rmk: complete series}
    Let $\D=|E|$, then $D_j'(P)=L(E-jP)$. Therefore $j$ is an $(|E|,P)$-order if and only if $\ell(E-jP)=\ell(E-(j+1)P)+1$. Applying the Riemann-Roch theorem, we see that this is also equivalent to $L(C-E+(j+1)P)=L(C-E+jP)$, for $C$ a canonical divisor on $\X$. 
\end{remark}

Now, let us see an application of Remark \ref{rmk: complete series}.

\begin{example}
    Let $g$ be the genus of $\X$ and $\D:=|E|$ with $d=\deg(E)\geq2$. We will show how by applying  Remark \ref{rmk: complete series} we may easily compute some $(\D,P)$-orders for any point $P\in\X$. Consider $0\leq i\leq d-2g$, then $\deg(C-E+(i+1)P)<0$ and, as a consequence, $\{0\}=L(C-E+(i+1)P)\supseteq L(C-E+iP)$ and we may deduce that these $i$'s are $(\D,P)$-orders. In particular, $\D$ is base-point-free and $j_i=i$ for $0\leq i\leq d-2g$.
\end{example}

In the following example, we will see how any given sequence of non-negative integers can be made into a sequence of Hermitian invariants for a linear series on some curve and some point on it.

\begin{example}
    Let $\ell_0<\dots<\ell_r$ be non-negative integers. We will exhibit a curve $\Y$, a point $P_0\in\Y$ and a linear series $\Fs$ on $\Y$ such that the $\ell_i$'s for $i=0,\dots,r$ are exactly the $(\Fs,P_0)$-orders. We set $\Y:=\P^1(\F)$ and $x$ a transcendental element over $\F$, so that we can identify $\F(\Y)$ with $\F(x)$. Moreover, $P_{\infty}:=(0:1)$, and $P_a:=(1:a)$ for $a\in\F$. We also assume $\mathrm{div}(x)=P_0-P_{\infty}$, and $v_{P_a}(x-a)=1$ for $a\in\F$. We define $E:=\ell_rP_\infty$ and $\Fs':=\Span\left(x^{\ell_0},\dots,x^{\ell_r}\right)\subseteq\F(x)$. Clearly, $\Fs'$ is a linear subspace of $\mathcal{L}(E)$, since $E+\div(x^{\ell_i})=\ell_iP_0+(\ell_r-\ell_i)P_\infty$ and $\ell_r-\ell_i\geq0$ for $i=0,\dots,r$. In particular, $\Fs:=\{E+\div(f)\mid f\in\Fs'\}$ is a $g_{\ell_r}^r$ linear series on $\Y$. Moreover, we see that $\Fs_{\ell_i}'=\Fs'\cap L(E-\ell_iP_0)=\Span\left(x^{\ell_i},\dots,x^{\ell_r}\right)$ has dimension $r-i$ for all $i=0,\dots,r$ and thus $\ell_0,\dots,\ell_r$ are $(\Fs,P_0)$-orders. Since we have exhibited $r+1$ of them, they are precisely the orders. Moreover, we see that $b(P)=0$ for every $P\neq P_0$, thus the base locus of $\Fs$ is $B=\ell_0P_0$. We may also see that $\Fs$ is complete if and only if $\ell_0=0$ (that is $\Fs$ is base-point-free) and $\ell_r=r$. In this case, the image $\phi(\Y)$ of the curve $\Y$ via the morphism $\phi=(x^{\ell_0}:\dots:x^{\ell_r})$ associated to $\Fs$ (see the beginning of the next section), is the so-called \textit{rational normal curve} in $\P^r$.
\end{example}

We are now ready to introduce the \textit{Weierstrass semigroup}. Notice that, by Remark \ref{rmk: complete series}, $j$ is an $(|E|,P)$-order if and only if there is \textbf{no} $f$ in $L(C-E+(j+1)P)$ such that $v_P(C-E)+v_P(f)=-(j+1)$, where $C$ is a canonical divisor on the curve $\X$. This motivates the following definition.

\begin{definition}
    Let $D$ be a divisor on $\X$. We say that a non-negative integer $\ell$ is a \textit{$(D,P)$-gap} if there is no $f$ in $L(D+\ell P)$ such that $v_P(D)+v_P(f)=-\ell$.
\end{definition}

This definition, together with the above observation, tells us that $\ell$ is a $(D,P)$-gap if and only if $\ell-1$ is a $(|C-D|,P)$-order. Notice that in particular $(D,P)$-gaps are always $\geq1$. 

\begin{definition}
    We call \textit{Weierstrass gap} at $P$ any $(0,P)$-gap and we denote by $G(P)$ the set of Weierstrass gaps. The \textit{Weierstrass semigroup} at $P$ is the set $H(P):=\N\backslash G(P)$ and its elements are called \textit{Weierstrass non-gaps} at $P$.
\end{definition}

The following are direct consequences. Firstly, we have that $\#G(P)=g$, and this is known as \textit{Weierstrass gap theorem}. Indeed, $\#G(P)=\#\{(|C|,P)\mbox{-orders}\}=\dim(|C|)+1=\ell(C)=g$. Moreover, we have that $h$ is a Weierstrass non-gap ($h\in H(P)$) if and only if there exists $f_h$ in $L(hP)$ such that $v_P(f_h)=-h$, or equivalently if there exists $f_h$ in $\F(\X)$ such that $\div_\infty(f_h)=hP$.
As a consequence of this last equivalence, we are also reassured that $H(P)\subseteq\N$ is indeed a semigroup, since we may see that it is closed under the sum. Notice finally that $H(P)$ contains every integer $h\geq 2g$, where $g$ is the genus of $\X$. Indeed recall that $\ell$ is a Weierstrass gap if an only if $\ell-1$ is a $(|C|,P)$-order and thus $\ell-1\leq\deg(|C|)=\deg(C)=2g-2$, therefore every $h\geq 2g$ has to be a Weierstrass non-gap.

Let us write the elements of the Weierstrass semigroup $H(P)$ as a strictly increasing sequence $n_i(P)$ for $i\geq0$ (we will also just write $n_i$ for $n_i(P)$). Then, in particular, $n_0(P)=0$ and thus one can inductively see that $\ell(n_i(P)P)=i+1$ (for this, one uses the facts that $\ell(0)=1$, and $h\in H(P)$ if and only if $\ell(hP)=\ell((h-1)P)+1$, which is proved with the same argument as in Remark \ref{rmk: complete series}). Notice that $\{0,\dots,2g-1\}=G(P)\cup\{n_0,\dots,n_{g-1}\}$ and therefore $n_i(P)=g+i$ for $i\geq g$.

\begin{remark}
    If $\X$ is a curve of genus $g\geq1$, then $|C|$ is base-point free. Indeed, the element $1$ must be a gap for every $P\in\X$ (otherwise we would have that $H(P)=\N$ and thus $0=\#G(P)=g$), meaning that $0$ is a $(|C|,P)$-order for every $P$.
\end{remark}

The following is an example of linear series coming from the Weierstrass semigroup. This kind of linear series have important applications in the literature on maximal curves over finite fields, see \cref{rem:appl}.

\begin{example}
\label{ex:nr}
 Let $\X$ be a curve, $P$ be a point on $\X$ and $n_i:=n_i(P)$ be the elements of $H(P)$. Consider the complete linear series $\D:=|n_rP|$, which is by definition a $g_{n_r}^r$ on $\X$. Let us see that $\D$ is base-point-free. Clearly $P$ cannot be a base point, indeed since $n_r\in H(P)$, we find $f\in L(n_rP)$ such that $v_P(f)=-n_r$ and thus $b(P)=0$. The same holds for $Q\neq P$ since $D:=n_rP+\div(1)\in \D$ and $v_Q(D)=0$. Now let us see that the $(\D,P)$-orders are exactly the integers $n_r-n_i$ with $i=0,\dots,r$. We thus need to show that $\D_{n_r-n_i}(P)\supsetneq\D_{n_r-n_i+1}(P)$. Since $n_i\in H(P)$ we find $f_i\in L(n_iP)$ such that $\div(f_i)=\div_0(f_i)-n_iP$. Then $n_rP+\div(f_i)=(n_r-n_i)P+\div_0(f_i)$ is indeed an element of $\D_{n_r-n_i}(P)$ which is not in $\D_{n_r-n_i+1}(P)$
\end{example}

\begin{remark}
\label{rem:appl}
    Linear series arising from Weierstrass non-gaps, as those discussed in \cref{ex:nr}, played a key role in the proofs of uniqueness, up to $\F_{q^2}$-isomorphism, of maximal curves having the largest possible genus (that is, the Hermitian curve) \cite{RS}, and the second-largest possible genus \cite{FGT}, \cite{AT}. The problem of characterizing $\F_{q^2}$-maximal curves with the third-largest possible genus is still open. 
\end{remark}

\medskip

\subsection{Morphisms and linear series}

We start with some recalling. For $P\in\X$ a point, we say that a function $f\in\F(\X)$ is \textit{regular at $P$} if $v_P(f)\geq0$. When this is the case, then there is a unique $a_f\in\F$ such that $v_P(f-a_f)\geq0$ and one sets $f(P):=a_f$. Any morphism $\phi\colon\X\rightarrow\P^r$ may be written as $\phi(P)=\left((t^{e_P}f_0)(P):\dots:(t^{e_P}f_r)(P)\right)$ where $f_0,\dots, f_r\in\F(\X)$ are not all zero, $t$ is a \textit{local parameter} at $P\in\X$ (that is a rational function $t\in\F(\X)$ such that $v_P(t)=1$) and $e_P:=-\min\left\{v_P(f_0),\dots,v_P(f_r)\right\}$. We then notice that the functions $t^{e_P}f_i$ are by definition regular at $P$ and therefore $\phi$ is well defined. The rational functions $f_0,\dots,f_r$ are called (homogeneous) \textit{coordinates} (they are uniquely determined by $\phi$ up to a factor in $\F(\X)^*$) of $\phi$ and one writes $\phi=(f_0:\dots:f_r)$.

Following \cite{torres2000}, we will show how to associate a morphism to a linear series and viceversa. Let $\D\cong\P(\D')$ be an $r$-dimensional linear series on $\X$.

\medskip

\textbf{Associating a morphism to $\D$.} Recall that $\dim(\D_{b(P)+1})=\dim(\D)-1$, that is $\D_{b(P)+1}\subsetneq\D$ is a hyperplane and can thus be seen as a point in the dual $\D^*\cong\P(\D')^*$. We then have the map $$\phi_\D\colon\X\rightarrow\D^*\cong\P(\D')^*,\quad P\mapsto \D_{b(P)+1}.$$ Let us give a more concrete description of this morphism in terms of its homogeneous coordinates. Let $\{f_0,\dots,f_r\}$ be an $\F$-basis of $\D'$ and $t$ be a \textit{local parameter} at $P\in\X$. For $f\in\D'\backslash\{0\}$, we see that $E+\div(f)$ belongs to $\D_{b(P)+1}$ if and only if $v_P(t^{v_P(E)-b(P)}f)\geq1$ which is in turn equivalent to having that $\left(t^{v_P(E)-b(P)}f\right)(P)=0$. Writing $f=\sum_{i=0}^ra_if_i$ with $(a_0:\dots:a_r)\in\P^r$, we have $$\D_{b(P)+1}\cong\left\{(a_0:\dots:a_r)\in\P^r\mid\sum_{i=0}^r\left(t^{v_P(E)-b(P)}f_i\right)(P)a_i=0\right\}$$$$\cong\left(\left(t^{v_P(E)-b(P)}f_0\right)(P):\dots:\left(t^{v_P(E)-b(P)}f_r\right)(P)\right)\in\P^r.$$ We then have that the morphism $\phi_{f_0,\dots,f_r}:=(f_0:\dots:f_r)$ gives a coordinate description of $\phi_\D$. Moreover, for any other $\F$-basis $\{g_,\dots,g_r\}$ of $\D'$, there exists $T\in\mathrm{Aut}(\P^r)$ such that $\phi_{g_0,\dots,g_r}=T\circ\phi_{f_0,\dots,f_r}$. Therefore, up to projective equivalence, $\phi_{f_0,\dots,f_r}$ is uniquely determined by $\D$ and we call it \textit{morphism associated to $\D$}.

\medskip

\textbf{Linear series from morphisms.} In this part, we explain how a morphism from a curve $\X$ to a projective space defines a corresponding linear series.

Let $\phi=(f_{0}:\dots : f_{r}): \X \rightarrow \P^r$ be a morphism on $\X$. Let $e_{P}:=-\min\{v_{P}(f_{0}),\dots,v_{P}(f_{r})\}$ where $f_{0},\dots,f_{r}\in \F(\X)$ are all non-zero rational functions, $P\in \X$ and $v_{P}(f)\ge 0.$ Then, $e_{P}\neq 0$ for finitely many $P\in \X$ and there exist a divisor $E=E_{f_0,\dots,f_r}$ such that $v_{P}(E):=e_{P}$. This makes it possible to define a linear series generated by the coordinate functions of the morphism. Let
\[
\D'=\langle f_0,\dots,f_r\rangle \subseteq \mathcal{L}(E).
\]
The associated linear series is
\[
\D_{f_0,\dots,f_r}=\{E+\operatorname{div}(f)\mid f\in \D'\setminus\{0\}\}\subseteq |E|.
\]

A key point is that this linear series is base-point-free. Indeed, at each point \(P\), one may choose a coordinate \(f_{i_0}\) for which \(e_P=-v_P(f_{i_0})\), and then \(E+\div(f_{i_0})\) has multiplicity zero at $P$. Thus the construction yields a canonically attached base-point-free linear series. Moreover, if we replace ${f_0,\dots,f_r}$ by a projectively equivalent morphism, $\D_{f_0,\dots,f_r}$ does not change. In addition, multiplying all coordinates by the same non-zero rational function also leaves the resulting linear series unchanged. In addition, the construction is invariant under projective equivalence: if $\phi_1=T\circ \phi$ for some $T\in \operatorname{Aut}(\P^r)$, then $\D_{\phi_1}=\D_\phi$. Likewise, multiplying all homogeneous coordinates by a common non-zero rational function does not change the resulting linear series. Consequently, $\D_\phi$ depends only on the projective equivalence class of the morphism.

From now on, we will use the linear series $\D_{f_0,\dots f_r}$ briefly as $\D_{\phi}$ for a morphism $\phi=(f_0:\dots:f_r)$. 

When the morphism is \textit{non-degenerate}, meaning that its image is not contained in any hyperplane, the coordinate functions $f_0,\dots,f_r$ are linearly independent, and therefore the associated linear series has dimension $r$. In this case, the section gives a very geometric description of $\D_\phi$: it is exactly the family of pull-backs of hyperplanes in $\P^r$,
\[
\D_\phi=\{\phi^*(H)\mid H \text{ is a hyperplane in } \P^r\}.
\]

Thus, a linear series attached to a morphism may be interpreted as the hyperplane section series induced on the curve by the map $\phi$. This is one of the main ideas of the section. We can give one of the several basic properties of these pull-backs as the degree of the associated linear series is shown to satisfy \[
\deg(\D_\phi)=\deg(\phi)\deg(\phi(\X)) \text{ where } \phi:\X\rightarrow \P^r \text{ be an non-degenerate morphism and birational}.
\]
Then, this formula identifies the degree of the induced linear series as the product of the degree of the morphism and the degree of the image curve. In particular, if $\phi$ is birational, then $(\deg(\D_\phi)=\deg(\phi(\X))$; if $\X\subseteq \P^r$ and $\phi$ is the inclusion, then $\D_\phi$ coincides with the classical linear series of hyperplane sections $\X\cdot H$. Thus the general construction recovers, as a special case, the standard projective linear series attached to an embedded curve. In other words, $
\D_\phi=\{\X\cdot H \mid H \text{ hyperplane}\}$, this construction generalizes the familiar notion of hyperplane sections of an embedded curve.

\begin{example}\label{morphtolinser}
We illustrate the construction of a linear series from a morphism in the simplest non-trivial case. Let $\X=\P^1$
and let $t$ be the standard affine coordinate on $\X$. Denote by $P_\infty$ the point at infinity. Then
\[
\div(t)=P_0-P_\infty,
\]
so, $t$ has a pole of order $1$ at $P_\infty$, while $t^2$ has a pole of order $2$ there.
Consider the morphism
\[
\phi=(1:t:t^2)\colon\X\longrightarrow \P^2.
\]
Following the general construction, for each $P\in \X$ we define
\[
e_P=-\min\{v_P(1),v_P(t),v_P(t^2)\}.
\]
If $P\neq P_\infty$, then $1$, $t$, and $t^2$ are regular at $P$, hence $
e_P=0$.
At $P_\infty$, we have $
v_{P_\infty}(1)=0,\quad v_{P_\infty}(t)=-1,\quad v_{P_\infty}(t^2)=-2$,
and therefore $
e_{P_\infty}=2$.
It follows that the divisor $E$ associated with $\phi$ is $E=2P_\infty$.
The corresponding vector space is
\[
\D'=\langle 1,t,t^2\rangle \subseteq L(2P_\infty),
\]
and the associated linear series is
\[
\D_\phi=\bigl\{\,2P_\infty+\div(a+bt+ct^2)\mid(a:b:c)\in \P^2\,\bigr\}.
\]
Thus $\D_\phi$ is a linear series of degree $2$ and projective dimension $2$, i.e. a $g^2_2$. We verify that $\D_\phi$ is base-point-free. If $P\neq P_\infty$, then the divisor $
2P_\infty+\div(1)=2P_\infty$
does not contain $P$. On the other hand, $
\div(t^2)=2P_0-2P_\infty$, so, 
$2P_\infty+\div(t^2)=2P_0,$
which does not contain $P_\infty$. Hence no point of $\X$ is contained in every divisor of the series, and then $\D_\phi$ is base-point-free. Then, we determine the image of $\phi$. Writing homogeneous coordinates on $\P^2$ as $(X_0:X_1:X_2)$, the image satisfies $X_1^2=X_0X_2,$
since $( X_0:X_1:X_2 )=(1:t:t^2).$
Moreover, $\phi$ is birational onto its image, because $t=X_1/X_0$ on the affine chart $X_0\neq 0$. Hence $
\deg(\phi)=1.$ The general degree formula gives $
\deg(\D_\phi)=\deg(\phi)\deg(\phi(\X))=1\cdot 2=2.$
Finally, let
\[
H:\ aX_0+bX_1+cX_2=0
\]
be a hyperplane in $\P^2$. Then its pull-back under $\phi$ is $
\phi^*(H)=2P_\infty+\div(a+bt+ct^2)$.
For instance, if we take
\[
H:\ X_2-X_0=0,
\]
then $
a+bt+ct^2=t^2-1=(t-1)(t+1)$,
and therefore 
$\phi^*(H)=P_1+P_{-1}.$
This explicitly exhibits the associated linear series as the hyperplane-section series attached to the morphism:
\[
\D_\phi=\{\phi^*(H)\mid H \text{ a hyperplane in } \P^2\}.
\]
In other words, $P_1+P_{-1}\in \D_{\phi}$ where $H: X_2-X_0=0.$
\end{example}
This example shows how a morphism $\X\to \P^r$ gives rise to a base-point-free linear series on $\X$, and how the geometric properties of the image are reflected in the divisor-theoretic structure of the curve. 

\medskip

\textbf{Relation between linear series and morphisms.} Fix an integer $r\geq 1$. Following Torres, let
\[
\sL=\sL_r:=\{\D^{B} : \D \text{ is a linear series with } \dim(\D)=r\},
\]
where $\D^{B}$ denotes the base-point-free linear series obtained from $\D$ by subtracting its base locus, and let
\[
\M=\M_r:=\{\langle \phi \rangle : \phi:\X\to \P^r \text{ is a non-degenerate morphism}\},
\]
where $
\langle \phi \rangle=\{T\circ \phi : T\in \operatorname{Aut}(\P^r)\} $
is the projective equivalence class of $\phi$. Thus $\sL_{r}$ is the set of base-point-free linear series of projective dimension $r$, while $\M_r$ is the set of projective equivalence classes of non-degenerate morphisms from $\X$ to $\P^r$.

Considering previous parts, we can define two natural maps as follows. First map sends a base-point-free linear series to the projective equivalence class of its associated morphism:
\[
M=M_r:L\longrightarrow M,
\qquad
\D^{B}\longmapsto \langle \phi_{\D^{B}}\rangle.
\]
Here $\phi_{\D^{B}}$ is the non-degenerate morphism determined, up to an automorphism of $\P^r$, by a basis of the vector space defining $\D^{B}$. The second map sends a non-degenerate morphism to its associated linear series:
\[
L=L_r:\M\longrightarrow \sL,
\qquad
\langle \phi\rangle \longmapsto \D_{\phi}.
\]
The key observation is that these two maps are inverse to one another. More precisely, Torres shows that $M\circ L=\mathrm{id}_{\M}$
by definition, and $L\circ M=\mathrm{id}_{\sL}.$

It follows that base-point-free linear series of dimension $r$ are in one-to-one correspondence with projective equivalence classes of non-degenerate morphisms $\X\to \P^r$.

\begin{remark}\label{linser:hyper}
If $\D$ is a linear series of dimension $r$, and $\phi:\X\to \P^r$ is the non-degenerate morphism determined by a basis of the vector space defining $\D$, then the associated base-point-free series $\D^{B}$ is exactly the family of pull-backs of hyperplanes:
\[
\D^{B}=\{\phi^{*}(H): H \text{ hyperplane in } \P^r\}\subseteq |E-B|.
\]
Thus the divisors of the base-point-free series are precisely the hyperplane sections induced by the corresponding morphism. 
\end{remark}

For the illustration of this process, we can give an example as follows. 
\begin{example} Let $\X=\P^1$ and $t$ be a standard affine coordinate on $\P^1$. Consider the linear series $\D=|2P_{\infty}|$. If we consider Example \ref{morphtolinser}, we can say that $L(2P_{\infty})=\langle 1,t,t^2\rangle.$ Then, the series $\D$ is base-point-free $g^2_2.$ The associated morphism is $\phi=(1:t:t^2):\P^1\longrightarrow \P^2$ and its image is the smooth conic $X_1^2=X_0X_2.$ Now let $H:\ aX_0+bX_1+cX_2=0$
be a hyperplane in $\P^2$. Then the pull-back of $H$ under $\phi$ is
\[
\phi^{*}(H)=2P_{\infty}+\div(a+bt+ct^2).
\]
As $(a:b:c)$ in $\P^2$, these divisors are precisely the elements of the linear series $\D$. Therefore
\[
\D=\{\phi^{*}(H): H \text{ a hyperplane in } \P^2\},
\]
which is exactly the content of Remark \ref{linser:hyper} in this concrete case. For instance, if one takes the line $
H:\ X_2-X_0=0$,
then $
a+bt+ct^2=t^2-1=(t-1)(t+1)$,
and hence
\[
\phi^{*}(H)=2P_{\infty}+\operatorname{div}(t^2-1)=P_1+P_{-1}.
\]
Thus the divisor $P_1+P_{-1}$ is one element of the linear series $\D=|2P_{\infty}|$, obtained as the pull-back of a hyperplane section of the image conic. This makes explicit the general principle that a base-point-free linear series may be interpreted as the hyperplane-section series induced by its associated morphism.

In other words, starting from the base-point-free linear series $|2P_{\infty}|$, one obtains the associated morphism $\phi=(1:t:t^2):\P^1\to \P^2$.
Conversely, applying the construction given in the part \textit{Linear Series from Morphism} to $\phi$, one recovers the same linear series:
$\D_{\phi}=|2P_{\infty}|$. Thus this example realizes concretely the correspondence between base-point-free linear series of dimension $r$ and projective equivalence classes of non-degenerate morphisms to $\P^r$.
\end{example}

\section{The order bound for duals of one-point AG codes}
\label{sec:orderbound}

\subsection{Recalls on AG codes and their duals}

Consider a projective, irreducible and non-singular algebraic curve $\mathcal{C}$ of genus $g$ defined over a finite field $\F_q$. Let $G$ and $D$ be two divisors whose supports are disjoint and such that $D$ is a sum of $n$ distinct points of $\mathcal{C}(\F_q)$: $D = \sum_{i=1}^n P_n$. Recall that we denote by $\mathcal{L}(G)$ the Riemann--Roch space associated to the divisor $G$, \emph{i.e.} $\mathcal{L}(G) = \{ f \in \F_q(\mathcal{C})^* \mid \div(f) \geq -G \} \cup \{0\}$.

By the evaluation map
\begin{equation*}
	        \begin{split}
	            ev_D\colon \mathcal{L}(G) & \longrightarrow \mathbb{F}_q^n\\
	            f & \longmapsto (f(P_1),\dots,f(P_n)), 
	        \end{split}
\end{equation*}
we define the AG\footnote{\emph{Algebraic Geometry}} code $C_L(D,G)$ on $\mathcal{C}$ as the image
\[
C_L(D,G) := \{ ev_D(f) \, | \, f \in \mathcal{L}(G) \} . 
\]
We say that $C_L(D,G)$ is a \emph{one-point AG code} when the support of $G$ is reduced to a single point, \emph{i.e.} $G = hP$ where $h$ is an integer and $P$ is a point.

Recall that we define the (Hamming) distance between two elements $c=(c_1, \dots, c_n)$ and $c'=(c'_1, \dots, c'_n) \in\F_q^n$ to be
\[
d(c,c') := \# \{ 1 \leq i \leq n \mid c_i \neq c'_i\} \, . 
\]
Note that $d(c,c') = d(c-c',0)$ and we define the \emph{minimum (Hamming) distance} of a linear code $C \in \F_q^n$ to be the integer
\[
d(C) := \min \{ d(c,0) \mid c \in C \setminus \{ 0 \} \} \geq 1 \, .
\]

In general, we can compute a lower bound for the minimum distance of the AG code $C_L(D,G)$ that depends on the parametrizing divisors.

\begin{definition}[Goppa Bound {\cite[Definition 2.2.4]{stichtenoth2009algebraic}}]
    Let $C_L(D,G)$ be an AG code parametrized by a divisor $G$. The integer $d_* := n - \deg(G) = \deg(D-G)$ is called the \emph{Goppa bound} of $C_L(D,G)$.
\end{definition}
Suppose $\deg(G) < n$. Since $L(G-D) =0$, the evaluation map $ev_D\colon \mathcal{L}(G) \mapsto \mathbb{F}^n_q$ is injective and the dimension $k$ of the AG code $C_L(D,G)$ verifies
\[
k = \dim_{\F_q}(\mathcal{L}(G)) \geq \deg(G) -g + 1 \, .
\]
In any case, its minimum distance $d$ is bounded from below by the Goppa Bound:
\[
d \geq n-\deg(G) \, .
\]
Let $c = (c_1, \dots, c_n)$ and $c'=(c'_1, \dots, c'_n)$ be two elements of $\F_q^n$. We define their canonical scalar product on $\F_q^n$ as the sum
\[
c \cdot c' := \sum_{i=1}^n c_ic'_i \, .
\]
The \emph{dual} of a linear code $C \subset \F_q^n$, denoted $C^\perp$, is the orthogonal space of $C$ for this scalar product.
 Through a combination of Serre's duality theorem \cite[Theorem 1.5.14]{stichtenoth2009algebraic} and the residue theorem \cite[Corollary 4.3.3]{stichtenoth2009algebraic}, we can express the dual of $C_L(D,G)$ as an AG code on the curve $\mathcal{C}$.

\begin{theorem}[{\cite[Proposition 2.2.10]{stichtenoth2009algebraic}}]\label{thm:AGdual}
    Let $\omega$ be a differential form on $\mathcal{C}$ whose valuation and residue at every point $P \in \Supp(D)$ is $v_P(\omega) = -1$ and $res_P(\omega) = 1$ respectively. Then
    \[
    C_L(D,G)^\perp = C_L(D,D-G+\div(\omega)).
    \]
\end{theorem}
Note that by the approximation Theorem \cite[Theorem 1.6.5]{stichtenoth2009algebraic}, such a divisor $\omega$ always exists and it is easy to check that the supports of $D$ and $D-G + \div(\omega)$ remain disjoint.

Let $C_L(D,G)$ be an AG code and consider its dual $C_L(D,G)^\perp = C_L(D,G')$ where $G'$ is a divisor of the form $G'=D-G+\div(\omega)$. By properties of the dual, the dimension of $C_L(D,G)^\perp$ is $\dim_{\F_q}( C_L(D,G)^\perp)= n- \dim_{\F_q}(C_L(D,G)^\perp$ and if the evaluation map $ev_D\colon L(D-G+\div(\omega)) \mapsto \mathbb{F}^n_q$ is injective then 
\[
\dim_{\F_q}(C_L(D,G)^\perp) \geq n - \deg(G) +g -1 \, .
\]
Applying the Goppa bound to $C_L(D,D-G+\div(\omega))$ gives a lower bound for the minimum distance $d(C_L(D,G)^\perp)$ of $C_L(D,G)^\perp$:
\[
d(C_L(D,G)^\perp) \geq \deg(G) - 2g + 2 \, .
\]
Before moving on to the Feng--Rao bound, we will give two fundamental examples of one-point AG codes on algebraic curves, one on the projective line and the other on the Hermitian curve.

\begin{example}[{Reed--Solomon} code]
    Consider the projective line $\mathbb{P}^1$ on $\F_q$ whose rational points are $\mathbb{P}^1(\F_q) = \F_q \cup \{P_\infty\}$. Denote by $\{P_\alpha\}_{\alpha \in \F_q}$ the rational points of $\mathbb{P}^1(\F_q) \setminus \{P_\infty\}$ and consider the divisor of evaluation points $D := \sum_{\alpha \in \F_q} P_\alpha$. For every non-negative integer $k$, the \emph{Reed--Solomon code} of length $q$, parameter $k$ and on the alphabet of size $q$ is the one-point AG code on $\mathbb{P}^1_{\F_q}$ denoted by $\RS_q(k)$:
    \[
    \RS_q(k) := C_L(D,kP_\infty) \subset \F_q^q \, .
    \]
    Suppose $k < q$. Note that since the genus of $\mathbb{P}^1$ is zero, the dimension of the Riemann--Roch space $L(kP_\infty)$ is exactly $\dim_{\F_q}(L(kP_\infty)) = k + 1$ and thus
    \[
    \dim_{\F_q}(\RS_q(k)) = k + 1 \, .
    \]
    For the minimum distance $d$ of $\RS_q(k)$, the Goppa bound gives us a lower bound $d \geq q - k$ while the Singleton bound provides the upper bound $d \leq q + 1 - (\dim_{\F_q}(\RS_q(k))) = q-k$ and thus
    \[
    d = q-k \, .
    \]
    To compute the dual of $\RS_q(k)$, consider the differential form on $\mathbb{P}^1_{\F_q}$
    \[
    \omega :=  \frac{dx}{\prod_{\alpha \in \F_q} (x-\alpha)  }\,,
    \]
    with $P_\infty = \left \{ \frac{1}{x} = 0 \right \}$. Its associated canonical divisor is $\div(\omega) = (q-2) P_\infty - D$ and its residue at every point $P_\alpha \in \mathbb{P}^1(\F_q)$ is exactly $1$. By Theorem \ref{thm:AGdual}, this implies
    \[
    \RS_q(k)^\perp = \RS_q(q-2-k) \, .
    \]
\end{example}

\begin{example}[AG code on the Hermitian curve {\cite[Section 8.3]{stichtenoth2009algebraic}}]
    \label{ex:HermCode}
    Let $q$ be a power of a prime number, and consider the finite field $\F_{q^2}$. Let $\mathcal{H}$ denote the Hermitian curve, that is the plane curve on $\F_{q^2}$ of equation
    \[
    y^{q+1} = x^q + x \, .
    \]
    The number of $\F_{q^2}$-rational points of the Hermitian curve is exactly $q^3+1$, including the common pole of the functions $y$ and $x$ that we denote by $P_\infty$. By \cite[Lemma 6.4.4]{stichtenoth2009algebraic}, a canonical divisor is given by
    \[
    \div(dy)_{\mathcal{H}} = (q^2-q-2)P_\infty .
    \]
    Consider the set of affine rational points $\mathcal{H}(\F_q) \setminus \{ P_\infty \}$ and define its divisor of evaluation
    \[
    D := \sum_{P \in \mathcal{H}(\F_{q^2}) \setminus \{ P_\infty \}} P \, .
    \]
    Its degree is $\deg(D) = q^3$ and $D$ is linearly equivalent to $q^3 P_\infty$ through the following principal divisor:
    \[
    \div(y^{q^2}-y) = D - q^3P_\infty \, ,
    \]
    giving us the canonical divisor $\div\left(\frac{dy}{(y^{q^2}-y)} \right) = -D + (q^3 +q^2 - q -2)$. Note that at any point in the support of $D$, the residue of $\frac{dy}{(y^{q^2}-y)}$ is exactly $1$. Let $k$ be a non-negative integer, we now define the following one-point AG code with parameter $k$, on the alphabet of size $q^2$:
    \[
    C_\mathcal{H}(k) := C_L(D,kP_\infty) \, .
    \]
    The code $C_\mathcal{H}(k)$ has length $\deg(D) = q^3$ and its minimum distance $d(C_\mathcal{H}(k))$ is bounded by the Goppa bound
    \[
    d(C_\mathcal{H}(k)) \geq q^3 - k \, .
    \]
    If $k \geq q^3+q^2-q-2$, $C_\mathcal{H}(k)$ is the full code $C_\mathcal{H}(k) = \F_{q^2}^{q^3}$. Suppose that $0 \leq k \leq q^3+q^2-q -2$, then thanks to Theorem \ref{thm:AGdual} and the differential form $\frac{dy}{(y^{q^2}-y)}$, we have
    \[
    C_\mathcal{H}(k)^\perp = C_\mathcal{H}(q^3 + q^2 -q - 2 - k) \, .  
    \]
\end{example}

\subsection{Feng-Rao bound}

We will now see how, using the Weierstrass semigroup of a parametrizing point, the lower bound for the minimal distance of the dual of a one-point AG code may be improved. Define, as in the previous subsection, a projective, irreducible and non-singular algebraic curve $\mathcal{C}$ defined on a finite field $\F_q$, and choose a set of $n$ distinct rational points $P_1, \dots, P_n$ summed together in an evaluation divisor $D := \sum_{i=1}^n P_i$. Choose a rational point $P \in \mathcal{C}(\F_q) \setminus \{ P_1, \dots, P_n \}$ and consider the family of one-point AG codes $C_L(D,kP)$ on $\mathcal{C}$ defined for all integers $k$. One can easily see that, given two integers $k_1 \leq k_2$, it holds $C_L(D,k_1P) \subset C_L(D,k_2P)$ and $C_L(D,k_2P)^\perp \subset C_L(D,k_1P)^\perp$.

Consider now the Weierstrass semigroup at $P$
\[
H(P) = \{h_1:=0 < h_2 < \dots \} \subset \mathbb{N} \, ,
\]
and choose for every integer $i \geq 0$ a function $f_i \in \F_q(\mathcal{C})$ whose valuation at $P$ is exactly $\val_P(f_i)=i$ so that for all $i \geq 0$ we have
\[
L(h_iP) = \langle f_0,f_1, \dots,f_i \rangle \, . 
\]
For every integer $k \geq 0$, if we define $i_k := \max_i\{ h_i \leq k \}$, we see that $L(i_kP) = L(kP)$ which also implies $C_L(D,i_kP) = C_L(D,kP)$, and we can now consider the family of AG codes
\[
C_i := C_L(D,h_iP) \, ,
\]
for every integer $i \geq 0$. Our goal now is to find a lower bound for the minimum distance of $C_i^\perp$ that is better than the Goppa bound. Let $l \geq 0$ be an integer and suppose there exists $c \in C_l^\perp \setminus C_{l+1}^\perp$. We have the following scalar products
\begin{align*}
ev_D(f_l) \cdot  c &=0 \\
ev_D(f_{l+1}) \cdot c &\neq 0 \, .
\end{align*}
This implies that for every couple of non-negative integers $(i,j)$, we have
\begin{align*}
h_i + h_j \leq h_l &\implies ev_D(f_if_j) \cdot c = 0 \\
h_i + h_j = h_{l+1} &\implies ev_D(f_if_j) \cdot c \neq 0 \, ,
\end{align*}
and this observation motivates the following definition.

\begin{definition}
    Let $l \geq 0$ be an integer. Define
    \[
    N_l := \{ (i,j) \in \mathbb{N}^2 \mid h_i + h_j = h_{l+1} \} \, ,
    \]
    and $v_l := \#N_l$.
\end{definition}

Set $N$ an integer large enough so that $C_N = \F_q^n$ and chose $c \in C_l^\perp$ for any integer $l$. Consider the parity matrix of size $N \times n$

\[
H := \left( \begin{array}{c} ev_D(f_0) \\
    \hline
    
    \vdots \\
    \hline
    ev_D(f_N) \end{array} \right) \, .
\]
It is a full rank matrix and if we consider the diagonal matrix with $c$ written on its diagonal $D(c)$, then $rk(D(c)) = d(c,0)$ and the matrix $S(c) := HD(c)H^T$ also has $rk(S(c)) = d(c,0)$. Denote by $s_{i,j}(c)$ the $i$-th row and $j$-th column coordinate of $S(c)$. Then
\[
s_{i,j} = ev_D(f_if_j) \cdot c \, .
\]
It entails,
\begin{theorem}[\cite{HLP98}]
    \label{thm:FengRaoBound}
    Let $l$ be a positive integer and suppose $C_l \neq \F_q$ (\emph{i.e}, $C_l^\perp \neq \{0 \}$). Consider a codeword $c \in C_l^\perp \setminus C_{l+1}^\perp$, then
    \[
    d(c,0) \geq v_l \, .
    \]
\end{theorem}

\begin{definition}[\cite{FR94}]
    The \emph{Feng--Rao designed minimum distance} (or \emph{order bound}) of the AG code $C_l$ is defined as
    \[
    d_{ORD}(C_l) := \min \{ v_m \mid m \geq l \} \, .  
    \]
\end{definition}

By Theorem \ref{thm:FengRaoBound}, if $d(C_l^\perp)$ is the minimum distance of the dual $C_l$, then
\[
d(C_l^\perp) \geq d_{ORD}(C_l) \, .
\]
It can even be shown that this bound is an improvement from the Goppa Bound
\[
d_{ORD}(C_l) \geq \deg(h_lP) - 2g+2 \, ,
\]
where $g$ is the genus of $\mathcal{C}$.

\section{Conclusion} \label{section:conclusion}

In this paper, we have detailed Weierstrass semigroups and their connection to the order bound for related AG codes, highlighting the underlying geometry. 
In ongoing work, we are computing Weierstrass semigroups of the Garcia-Güneri-Stichtenoth curves $GGS_n$ \cite{GGS2010}, for $n\geq 5$. For $n=3$, that is for the Giulietti-Korchmáros curve, the Weierstrass semigroups are already known, see \cite{BMGKcurve2018}. For $n\geq 5$, the Weierstrass semigroups at all the $\F_{q^2}$-rational points of $GGS_n$ are known (see \cite{guneri2013automorphism}, \cite{bartoli2018ag}), but the problem of determining them for the remaining points is still open.
 Several central open problems remain in the study of Weierstrass semigroups of points on curves over finite fields. First, the classification of maximal and near‑maximal curves is not yet complete. Despite extensive work on the Hermitian curve and certain Galois subcovers, it is still unknown which numerical semigroups can occur as Weierstrass semigroups at rational or non‑rational points of general maximal curves over $\mathbb{F}_{q^2}$. Recent results show unexpectedly many distinct semigroup types even for well‑studied families, indicating that a full classification is highly non-trivial \cite{BMV26}, \cite{beelen2026weierstrass}. Second, there are interesting questions regarding realizability and arithmetic constraints. In particular, given a numerical semigroup, it is open to characterize when it can be realized as the Weierstrass semigroup of an $\mathbb{F}_q$-rational point on a curve defined over $\mathbb{F}_q$ \cite{zbMATH07562177}, \cite{MOYANOFERNANDEZ201946}. Third, work remains to be done in the optimization for algebraic‑geometric error‑correcting codes. While the Weierstrass semigroup at a rational point controls designed parameters and order bounds of one‑point and multi‑point AG codes, it remains open to systematically determine which semigroups (or families of semigroups) yield asymptotically or practically optimal codes over finite fields, and how semigroup invariants such as symmetry, telescopicity, or sparsity translate into improved decoding performance.

\section*{Acknowledgments}
This survey is part of the work of a research group, led by G.L. Matthews and L. Vicino, that started during the workshop \emph{Women in Numbers Europe 5 (WINE-5)}, held at the University of Split (Croatia), August 18-23, 2025. 

The workshop group is currently still working on a joint research project. The group would like to sincerely thank Maria Montanucci, who originally proposed the project and generously offered guidance and support during the initial steps of the work.

A. Barraud is supported by the grants ANR-21-CE39-0009-BARRACUDA and ANR-22-CPJ2-0047-01 from the French National Research Agency.

B. Gouthier is grateful to the Max Planck Institute for Mathematics in Bonn for their hospitality and financial support; she also thanks Heinrich-Heine-Universität Düsseldorf and Prof. Dr. Stefan Schröer for the support that permitted her attendance at the \textit{WINE-5} workshop.

G. L. Matthews is partially supported by NSF DMS-2502705 and the Commonwealth Cyber Initiative. 

L. Vicino is supported by a DFF-International Postdoc Grant (grant ID: 10.46540/5246-00028B) from Independent Research Fund Denmark.

\printbibliography

@article{zbMATH07562177,
 author = {{H. Charalambous and K. Karagiannis and S. Karanikolopoulos and A. Kontogeorgis}},
 title = {Weierstrass semigroups for maximal curves realizable as {Harbater}-{Katz}-{Gabber} covers},
 fjournal = {Advances in Geometry},
 journal = {Adv. Geom.},
 issn = {1615-715X},
 volume = {22},
 number = {3},
 pages = {445--450},
 year = {2022},
 doi = {10.1515/advgeom-2022-0014},
 zbMATH = {7562177},
 Zbl = {1493.14050}
}

@article{MOYANOFERNANDEZ201946,
title = {Generalized Weierstrass semigroups and their Poincaré series},
journal = {Finite Fields and Their Applications},
volume = {58},
pages = {46-69},
year = {2019},
issn = {1071-5797},
doi = {https://doi.org/10.1016/j.ffa.2019.03.005},
url = {https://www.sciencedirect.com/science/article/pii/S1071579719300267},
author = {J.J. Moyano-Fern\'{a}ndez and W. Ten\'{o}rio and F. Torres}
}

@misc{beelen2026weierstrasssemigroupsmaximalfermat,
      title={On Weierstrass semigroups of maximal Fermat function fields}, 
      author={P. Beelen and M. Montanucci and M. Frank vom Braucke},
      year={2026},
      eprint={2602.24015},
      archivePrefix={arXiv},
      primaryClass={math.AG},
      url={https://arxiv.org/abs/2602.24015}, 
}

@InProceedings{GL92,
author="Garcia, A.
and Lax, R. F.",
editor="Stichtenoth, H.
and Tsfasman, M.A.",
title="Goppa codes and Weierstrass gaps",
booktitle="Coding Theory and Algebraic Geometry",
year="1992",
publisher="Springer Berlin Heidelberg",
address="Berlin, Heidelberg",
pages="33--42",
isbn="978-3-540-47267-4"
}

@article {GGS2010,
    AUTHOR = {Garcia, A. and G\"{u}neri, C. and Stichtenoth, H.},
     TITLE = {A generalization of the {G}iulietti-{K}orchm\'{a}ros maximal
              curve},
   JOURNAL = {Adv. Geom.},
  FJOURNAL = {Advances in Geometry},
    VOLUME = {10},
      YEAR = {2010},
    NUMBER = {3},
     PAGES = {427--434},
      ISSN = {1615-715X},
   MRCLASS = {11G20 (14G15)},
  MRNUMBER = {2660419},
MRREVIEWER = {Michael E. Zieve},
       DOI = {10.1515/ADVGEOM.2010.020},
       URL = {https://doi.org/10.1515/ADVGEOM.2010.020},
}

@article {BMGKcurve2018,
    AUTHOR = {Beelen, P. and Montanucci, M.},
     TITLE = {Weierstrass semigroups on the {G}iulietti-{K}orchm\'{a}ros curve},
   JOURNAL = {Finite Fields Appl.},
  FJOURNAL = {Finite Fields and their Applications},
    VOLUME = {52},
      YEAR = {2018},
     PAGES = {10--29},
      ISSN = {1071-5797},
   MRCLASS = {11G20 (11R58 14H05 14H55)},
  MRNUMBER = {3807839},
MRREVIEWER = {Takehiro Hasegawa},
       DOI = {10.1016/j.ffa.2018.03.002},
       URL = {https://doi.org/10.1016/j.ffa.2018.03.002},
}

@book{stichtenoth2009algebraic,
  title={Algebraic function fields and codes},
  author={Stichtenoth, H.},
  year={2009},
  publisher={Springer}
}

@article{torres2000,
      title={The approach of Stöhr-Voloch to the Hasse-Weil bound with applications to optimal curves and plane arcs}, 
      journal={arXiv preprint},
      author={F. Torres},
      year={2000},
      eprint={math/0011091},
      archivePrefix={arXiv},
      primaryClass={math.AG},
      url={https://arxiv.org/abs/math/0011091}, 
}

@article{schmidt1939arithmetischen,
  title={Zur arithmetischen theorie der algebraischen funktionen. {II}. {A}llgemeine theorie der {W}eierstra{\ss}punkte},
  author={Schmidt, F.K.},
  journal={Mathematische Zeitschrift},
  volume={45},
  number={1},
  pages={75--96},
  year={1939},
  publisher={Springer}
}

@article{stohr1986weierstrass,
  title={Weierstrass points and curves over finite fields},
  author={St{\"o}hr, K.-O. and Voloch, J.F.},
  journal={Proceedings of the London Mathematical Society},
  volume={3},
  number={1},
  pages={1--19},
  year={1986},
  publisher={Oxford University Press}
}

@article{AT,
      author={Abd\'{o}n, M. and Torres, F.},
      title={On maximal curves in characteristic two},
      journal={Manuscripta Math.},
      volume={99},
      pages={39--53},
      date={1999}
}

@article{FGT,
      author={Fuhrmann, R. and Garcia, A. and Torres, F.},
      title={On maximal curves},
      journal={J. Number Theory},
      volume={67},
      pages={29--51},
      date={1997}
}

@article{GV,
  title={Weierstrass points on certain non-classical curves},
  author={Garcia, A. and Viana, P.},
  journal={Archiv der Mathematik},
  volume={46},
  pages={315--322},
  year={1986},
  publisher={Springer}
}

@article{RS,
      author={R\"{u}ck, H.G. and Stichtenoth, H.},
      title={A characterization of Hermitian function fields over finite fields},
      journal={J. Reine Angew. Math.},
      volume={457},
      pages={185--188},
      date={1994}
}

@article{FR94,
  title={A simple approach for construction of algebraic-geometric codes from affine plane curves},
  author={Feng, G.-L. and Rao, T.R.N.},
  journal={IEEE Transactions on Information Theory},
  volume={40},
  number={4},
  pages={1003--1012},
  year={1994},
  publisher={IEEE}
}

@article{BMV23,
title = {Weierstrass semigroups and automorphism group of a maximal curve with the third largest genus},
journal = {Finite Fields and Their Applications},
volume = {92},
pages = {102300},
year = {2023},
issn = {1071-5797},
author = {P. Beelen and M. Montanucci and L. Vicino}
}

@article{OrderBound,
      author={Beelen, P.},
      title={The order bound for general algebraic geometric codes},
      journal={Finite Fields and Their Applications},
      volume={13},
      pages={665-680},
      date={2007}
}

@article{LV22,
  title={Two-point AG codes from the Beelen-Montanucci maximal curve},
  author={Landi, L. and Vicino, L.},
  journal={Finite Fields and Their Applications},
  volume={80},
  pages={102009},
  year={2022},
  publisher={Elsevier}
}

@article{TwoPointCodes,
      author={{E. Barelli and  P. Beelen and M. Datta and V. Neiger and J. Rosenkilde}},
      title={Two-point codes for the generalized GK curve},
      journal={IEEE Transactions on Information Theory},
      volume={64},
      pages={6268-6276},
      date={2018}
    }

@article{LTV24,
  title={Two-point AG codes from one of the Skabelund maximal curves},
  author={Landi, L. and Timpanella, M. and Vicino, L.},
  journal={IEEE Transactions on Information Theory},
  year={2024}
}

@article{M04,
  title={Codes from the Suzuki function field},
  author={Matthews, G.L.},
  journal={IEEE Transactions on Information Theory},
  volume={50},
  number={12},
  pages={3298--3302},
  year={2004},
  publisher={IEEE}
}

@article{GKL93,
  title={Consecutive Weierstrass gaps and minimum distance of Goppa codes},
  author={Garcia, A. and Kim, S.J. and Lax, R.F.},
  journal={Journal of pure and applied algebra},
  volume={84},
  number={2},
  pages={199--207},
  year={1993},
  publisher={Elsevier}
}

@article{M01,
  title={Weierstrass pairs and minimum distance of Goppa codes},
  author={Matthews, G.L.},
  journal={Designs, Codes and Cryptography},
  volume={22},
  pages={107--121},
  year={2001},
  publisher={Springer}
}

@article{N2025,
title = {Non-isomorphic maximal function fields of genus $q-1$},
journal = {Finite Fields and Their Applications},
volume = {106},
pages = {102618},
year = {2025},
issn = {1071-5797},
author = {J. Niemann}
}

@misc{beelen2024some,
      title={Some families of non-isomorphic maximal function fields}, 
      author={{P. Beelen and M. Montanucci and J. Tilling Niemann and L. Quoos}},
      year={2024},
      eprint={2404.14179},
      archivePrefix={arXiv},
      primaryClass={math.NT},
      url={https://arxiv.org/abs/2404.14179}, 
}

@article{beelen2026weierstrass,
  title={Weierstrass semigroups and automorphism group of a maximal function field with the third largest possible genus, q$\equiv$ 1 (mod 3)},
  author={Beelen, P. and Montanucci, M. and Vicino, L.},
  journal={Finite fields and their applications},
  volume={109},
  pages={102701},
  year={2026},
  publisher={Elsevier}
}

@article{BMV26,
title = {Weierstrass semigroups and automorphism group of a maximal function field with the third largest possible genus, q$\equiv$ 0 (mod 3)},
journal = {Finite Fields and Their Applications},
volume = {110},
pages = {102729},
year = {2026},
issn = {1071-5797},
author = {P. Beelen and M. Montanucci and L. Vicino}
}

@article{BMZ21,
  author = "D. Bartoli and M. Montanucci and G. Zini",
  title = "{W}eierstrass semigroups at every point of the {S}uzuki curve",
  journal = "Acta Arithmetica",
  volume = "197",
  pages = "1-20",
  year = "2021"
}

@article{BM18GK,
  author = "P. Beelen and M. Montanucci",
  title = "{W}eierstrass semigroups on the {G}iulietti-{K}orchm\'{a}ros curve",
  journal = "Finite Fields and Their Applications",
  volume = "52",
  pages = "10-29",
  year = "2018"
}

@article{BLM21,
  author = "P. Beelen and L. Landi and M. Montanucci",
  title = "{W}eierstrass semigroups on the {S}kabelund maximal curve",
  journal = "Finite Fields and Their Applications",
  volume = "72",
  pages = "101811",
  year = "2021"
}

@article{G70,
  title={A new class of linear correcting codes},
  author="V. D. Goppa",
  journal={Problemy Peredachi Informatsii},
  volume={6},
  number={3},
  pages={24--30},
  year={1970},
  publisher={Russian Academy of Sciences, Branch of Informatics, Computer Equipment and~…}
}

@article{G71,
  title={A rational representation of codes and ({L},g)-codes},
  author="V. D. Goppa",
  journal={Problemy Peredachi Informatsii},
  volume={7},
  number={3},
  pages={41--49},
  year={1971},
  publisher={Russian Academy of Sciences, Branch of Informatics, Computer Equipment and~…}
}

@article{G77,
  title={Codes associated with divisors},
  author="V. D. Goppa",
  journal={Problemy Peredachi Informatsii},
  volume={13},
  number={1},
  pages={33--39},
  year={1977},
  publisher={Russian Academy of Sciences, Branch of Informatics, Computer Equipment and~…}
}

@article{G81,
  author = "V. D. Goppa",
  title = "Codes on algebraic curves",
  journal = "Doklady Akademii Nauk SSSR",
  volume = "259",
  pages = "1289-1290",
  year = "1981"
}

@article{G82,
  title={Algebraico-geometric codes},
  author="V. D. Goppa",
  journal={Izvestiya Rossiiskoi Akademii Nauk. Seriya Matematicheskaya},
  volume={46},
  number={4},
  pages={762--781},
  year={1982},
  publisher={Russian Academy of Sciences, Steklov Mathematical Institute of Russian~…}
}

@book{G88,
  author = "V. D. Goppa",
  title = "Geometry and codes",
  series = "Mathematics and its Applications (Soviet Series)",
  volume = "24",
  publisher = "Kluwer Academic Publishers Group",
  address = "Dordrecht",
  year = "1988"
}

@article{HLP98,
  title={Algebraic geometry codes},
  author={H{\o}holdt, T. and Van Lint, J.H. and Pellikaan, R.},
  journal={Handbook of Coding Theory},
  volume={1},
  number={Part 1},
  pages={871--961},
  year={1998},
  publisher={Elsevier Amsterdam}
}

@article{HK01,
  title={Goppa codes with {W}eierstrass pairs},
  author="M. Homma and S.J. Kim",
  journal={Journal of Pure and Applied Algebra},
  volume={162},
  number={2-3},
  pages={273--290},
  year={2001},
  publisher={Elsevier}
}

@article{HK05,
  title={Toward the determination of the minimum distance of two-point codes on a {H}ermitian curve},
  author="M. Homma and S.J. Kim",
  journal={Designs, Codes and Cryptography},
  volume={37},
  number={1},
  pages={111--132},
  year={2005},
  publisher={Springer}
}

@article{HK06,
  author = "M. Homma and S.J. Kim",
  title = "The complete determination of the minimum distance of two-point codes on a {H}ermitian curve",
  journal = "Designs, Codes and Cryptography",
  volume = "40",
  pages = "5-24",
  year = "2006"
}

@incollection{M04HE,
  author = "G.L. Matthews",
  title = "The {W}eierstrass semigroup of an {$m$}-tuple of collinear points on a {H}ermitian curve",
  booktitle = "Finite Fields and Applications",
  series = "Lecture Notes in Computer Science",
  publisher = "Springer-Verlag",
  address = "Berlin Heidelberg",
  volume = "2948",
  pages = "12-24",
  year = "2004"
}

@article{M05,
  author = "G.L. Matthews",
  title = "Weierstrass semigroups and codes from a quotient of the {H}ermitian curve",
  journal = "Designs, Codes, and Cryptography",
  volume = "37",
  pages = "473-492",
  year = "2005"
}

@article{bartoli2018ag,
  title={AG codes and AG quantum codes from the GGS curve},
  author={Bartoli, D. and Montanucci, M. and Zini, G.},
  journal={Designs, Codes and Cryptography},
  volume={86},
  number={10},
  pages={2315--2344},
  year={2018},
  publisher={Springer}
}

@article{montanucci2020ag,
  title={AG codes from the second generalization of the GK maximal curve},
  author={Montanucci, M. and Pallozzi Lavorante, V.},
  journal={Discrete Mathematics},
  volume={343},
  number={5},
  pages={111810},
  year={2020},
  publisher={Elsevier}
}

@article{s2025set,
  title={The set of pure gaps at several rational places in function fields},
  author={S. Castellanos, A. and Mendoza, E.A.R. and Tizziotti, G.},
  journal={Designs, Codes and Cryptography},
  volume={93},
  number={5},
  pages={1375--1400},
  year={2025},
  publisher={Springer}
}

@article{castellanos2024weierstrass,
  title={Weierstrass semigroups, pure gaps and codes on function fields},
  author={Castellanos, A.S. and Mendoza, E.A.R. and Quoos, L.},
  journal={Designs, Codes and Cryptography},
  volume={92},
  number={5},
  pages={1219--1242},
  year={2024},
  publisher={Springer}
}

@article{guneri2013automorphism,
  title={The automorphism group of the generalized Giulietti-Korchm{\'a}ros function field.},
  author={G{\"u}neri, C. and {\"O}zdemiry, M. and Stichtenoth, H.},
  journal={Advances in Geometry},
  volume={13},
  number={2},
  year={2013}
}

\end{document}